\documentclass[gmjn,parasubs,intlim,amssymb,righttag]{amsart} %
\usepackage{graphicx}
\begin{document}   

\newtheorem{thm}{Theorem}[section]   
\newtheorem{prop}[thm]{Proposition}
\newtheorem{lem}[thm]{Lemma}
\newtheorem{cor}[thm]{Corollary}

\newcommand{\R}{\mathbf{R}}


\title[Frequency control of SP forced Duffing's oscillator]{Frequency control of singularly perturbed forced Duffing's oscillator}


\author{Robert Vrabel} 
\address{Robert Vrabel, Institute of Applied Informatics, Automation and Mathematics, Faculty of Materials Science and Technology, Hajdoczyho 1,  917 01 Trnava, Slovakia}
\email{robert.vrabel@stuba.sk}

\author{Marcel Abas} 
\address{Marcel Abas, Institute of Applied Informatics, Automation and Mathematics, Faculty of Materials Science and Technology, Hajdoczyho 1,  917 01 Trnava, Slovakia}
\email{abas@stuba.sk} 

\date{{\bf January 10, 2011}(\bf JDCS Ref. 10-022, revised)}

\begin{abstract}
We analyze the dynamics of the forced singularly perturbed differential equation of Duffing's type. We explain the appearance of the large frequency nonlinear oscillations of the solutions. It is shown that the frequency can be controlled by a  small parameter at the highest derivative. We give some generalizations of results obtained recently by  B.S. Wu,  W.P. Sun and C.W. Lim, Analytical approximations to the double-well Duffing oscillator in large amplitude oscillations, \emph{ Journal of Sound and Vibration}, Volume 307, Issues 3-5, (2007), pp. 953-960.  The new method for an analysis of the nonlinear oscillations which is based on the dynamic change of coordinates is proposed.
\end{abstract}


\maketitle

 \section{Introduction}
Nonlinear oscillations comprise a large class of nonlinear dynamical systems, and arise naturally
from many physical systems such as mechanics,  chemistry, and engineering.
Also a variety of biological phenomena can be characterized as nonlinear oscillations, ranging from
heartbeat, neuronal activity, to population cycles (\cite{srebro}). 

The forced Duffing oscillator exhibits behavior, from limit cycles to chaos due to its nonlinear dynamics.
 When the periodic force that drives the system is large, chaotic behavior emerges and the phase space
 diagram is a strange attractor. In that case the behavior of the system is sensitive to the initial condition (\cite{Yeeal}).

In this work we focus our attention to the nonlinear oscillations in the context of the singularly perturbed 
forced oscillator of Duffing's type with a nonlinear restoring force
\begin{equation}\label{Duff_osc}
\epsilon^2\left(a^2(t)y'\right)'+f(y)=m(t),\quad 0<\epsilon<<1
\end{equation}
or rewriting to the autonomous system form 
{\setlength\arraycolsep{1pt}
\begin{eqnarray}
&\epsilon y'&=\frac wa\label{Duff_osc_system1}\\  
&\epsilon w'&=\frac{m(t)}{a}-\frac{f(y)}{a}-\epsilon\frac{a'}{a}w\label{Duff_osc_system2}\\
&t'&=1\label{Duff_osc_systemm3}
\end{eqnarray}}
where $a, m$ are the $C^1$ functions on the interval $\left\langle t_B,t_E\right\rangle,$ $a$ is positive and $f$ is a $C^1$ function on $\R.$

System (\ref{Duff_osc_system1}), (\ref{Duff_osc_system2}), (\ref{Duff_osc_systemm3}) is an example of a singularly perturbed system, because in the limit $\epsilon\rightarrow 0^+,$ it does not reduce to a differential equation of the same type, but to an algebraic-–differential reduced system
{\setlength\arraycolsep{1pt}
\begin{eqnarray*}
&0&=\frac wa\\  
&0&=\frac{m(t)}{a}-\frac{f(y)}{a}\\
&t'&=1.
\end{eqnarray*}}

Another way to study the singular limit $\epsilon\rightarrow 0^+$ is by introducing the new independent variable $\tau=\frac t\epsilon$ which transforms (\ref{Duff_osc_system1}), (\ref{Duff_osc_system2}), (\ref{Duff_osc_systemm3}) to the system
{\setlength\arraycolsep{1pt}
\begin{eqnarray*}
&\frac{\mathrm{d}y}{\mathrm{d}\tau}&=\frac wa\\  
&\frac{\mathrm{d}w}{\mathrm{d}\tau}&=\frac{m(t)}{a}-\frac{f(y)}{a}-\epsilon\frac{a'}{a}w\\
&\frac{\mathrm{d}t}{\mathrm{d}\tau}&=\epsilon.
\end{eqnarray*}}
Taking the limit $\epsilon\rightarrow 0^+,$ we obtain the so-called associated system (\cite{Jo})
{\setlength\arraycolsep{1pt}
\begin{eqnarray}
&\frac{\mathrm{d}y}{\mathrm{d}\tau}&=\frac wa\label{assoc_system1}\\  
&\frac{\mathrm{d}w}{\mathrm{d}\tau}&=\frac{m(t)}{a}-\frac{f(y)}{a}\label{assoc_system2}\\
&\frac{\mathrm{d}t}{\mathrm{d}\tau}&=0\label{assoc_system3}\quad \mathrm{i.e.}\quad t=t^*=const.
\end{eqnarray}}

The critical manifold $S$ is defined as a solution of the reduced system i.e. 
$$S:=\left\{(t,y,w):\quad t\in\left\langle t_B,t_E\right\rangle, f(y)=m(t),w=0\right\}$$
which corresponds to a set of equilibria for the associated  system (\ref{assoc_system1}), (\ref{assoc_system2}), (\ref{assoc_system3}).

We assume that
\begin{itemize} 
\item[$(A1)$] The critical manifold is S-shaped curve with two folds, i.e. it can be written in the form $t=\varphi(y),$ $t\in\left\langle t_B,t_E\right\rangle$ and the function $\varphi$
has precisely two critical points, one non-degenerate minimum $y_{\min}$ and one non-degenerate maximum $y_{\max}$ and let $y_{\min}<y_{\max}.$ Thus, the critical manifold can be broken up into three pieces $S_b,$ $S_m$ and $S_a,$ separated by the minimum and maximum (Fig. 1). These three pieces are defined as follows
{\setlength\arraycolsep{1pt}
\begin{eqnarray*}
&S_b&=\left\{(y,\varphi(y)):\quad y<y_{\min}\right\}\\
&S_m&=\left\{(y,\varphi(y)):\quad y_{\min}<y<y_{\max}\right\}\\
&S_a&=\left\{(y,\varphi(y)):\quad y_{\max}<y\right\}
\end{eqnarray*}}
\item[$(A2)$] $\varphi'(y)\neq0$ for $y\neq y_{\min},y_{\max}$

\item[$(A3)$] $\frac{\mathrm{d}f}{\mathrm{d}y}(y)<0$ for every $(t,y,0)\in S_m$ and $\frac{\mathrm{d}f}{\mathrm{d}y}(y)>0$ for every $(t,y,0)\in S_a\cup S_b.$

\end{itemize}

Let $t_{\min}=\varphi\left(y_{\min}\right),$ $t_{\max}=\varphi\left(y_{\max}\right).$  Denote by 
{\setlength\arraycolsep{1pt}
\begin{eqnarray*}
u_1(t)&=&\varphi^{-1}(t):\quad t\in\left\langle t_B,t_{\max}\right\rangle,\ y_{\max}\leq u_1(t) \\
u_2(t)&=&\varphi^{-1}(t):\quad t\in\left\langle t_{\min},t_{\max}\right\rangle,\ y_{\min}\leq u_2(t)\leq y_{\max}\\
u_3(t)&=&\varphi^{-1}(t):\quad t\in\left\langle t_{\min},t_E\right\rangle,\ u_3(t)\leq y_{\min}.
\end{eqnarray*}}
We divide the phase diagram of ({\ref{Duff_osc_system1}), (\ref{Duff_osc_system2}), (\ref{Duff_osc_systemm3}) into three charts, for $K_1,K_2,K_3,$ where
{\setlength\arraycolsep{1pt}
\begin{eqnarray*}
K_1&\subset&\left\langle t_B,t_{\min}\right)\\
K_2&\subset&\left( t_{\min},t_{\max}\right)\\
K_3&\subset&\left(t_{\max},t_E\right\rangle
\end{eqnarray*}}
are the compact sets.

\vglue3cm
\begin{figure}[ht]
\begin{center}
\includegraphics[scale=1, bb=6cm 17cm 15cm 20cm]{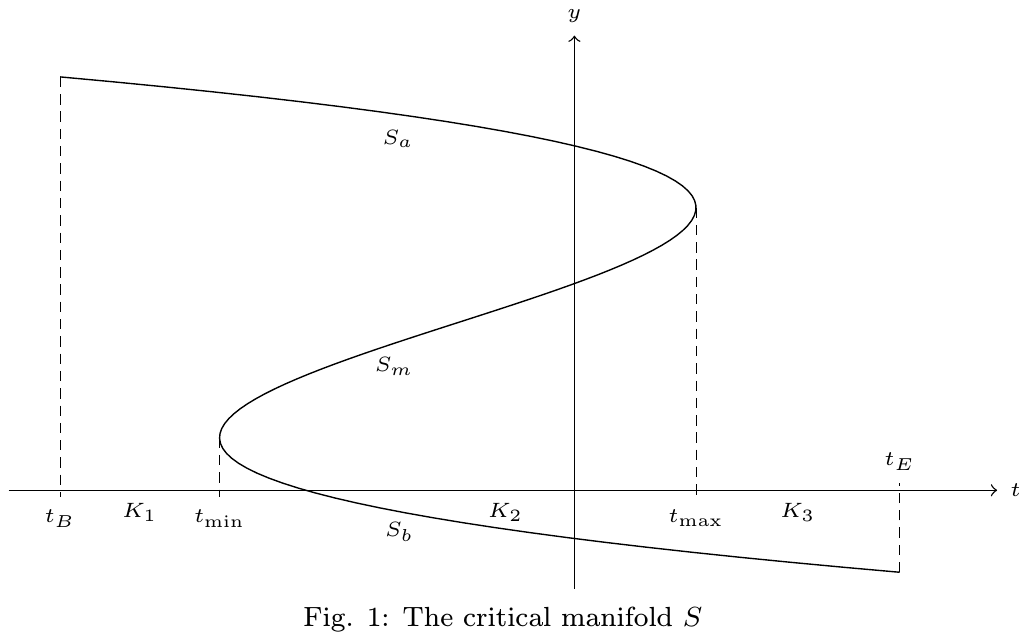}
\end{center}
\end{figure}

The situation considered here is principally different from the one in \cite{KrSz}, where two pieces of critical manifold are attracting and one is repelling.
In this paper,  two pieces $S_a$ and $S_b$ of the critical manifold are not normally hyperbolic (\cite{Jo})  and consequently the geometric singular perturbation theory developed by N. Fenichel (\cite{Fen}) is not applicable to our case. Indeed, all of the characteristic roots of associated system (\ref{assoc_system1}), (\ref{assoc_system2}), (\ref{assoc_system3}), $\lambda_{1,2}(t,y,w)=\pm a^{-1}(t)\sqrt{-\frac{\mathrm{d}f}{\mathrm{d}y}(y)},$ $\lambda_3=0,$ $(t,y,w)\in S_a\cup S_b$ lie on the imaginary axis. The piece $S_m$ is a normally hyperbolic manifold.

We generalize the results presented in \cite{WSL}, where unforced and undamped double-well Duffing oscillator with $\epsilon=1$ was examined. Moreover, the considerations below can be useful in the design of the high-frequency circuits (see e.g. \cite{HeHe, PChY}, and the references therein) and we introduce the parameter $\epsilon$ as a modeling tool for the frequency control of the oscillations.   

Our considerations relies on a suitable combination the phase-space analysis and the generalized polar coordinate transformations.

Consider the function
\begin{equation*}
H(t,y,w)=\frac12w^2+V(t,y),\quad V(t,y)=\int\limits_0^y f(s)\mathrm{d}s-m(t)y.
\end{equation*}

Let 
{\setlength\arraycolsep{1pt}
\begin{eqnarray*}
H^0(t) = \left\{ \begin{array}{ll}
V\left(t,u_1(t)\right) & \quad\textrm{for $t\in\left\langle t_B,t_{\min}\right)$}\\
V\left(t,u_2(t)\right) & \quad\textrm{for $t\in\left\langle t_{\min},t_{\max}\right\rangle$}\\
V\left(t,u_3(t)\right) & \quad\textrm{for $t\in\left(t_{\max},t_E\right\rangle$}. 
\end{array} \right.
\end{eqnarray*}}

We use the level surfaces $H(t,y,w)=H^\epsilon$  of $H$  with  
\begin{equation*}
H^\epsilon(t,y,w)=H^0(t)+\Delta(t)+ h^\epsilon(t,y,w) 
\end{equation*}
to characterize the trajectories of (\ref{Duff_osc_system1}), (\ref{Duff_osc_system2}), (\ref{Duff_osc_systemm3}), where  $h^\epsilon=O\left(\epsilon^\nu \right),$  $\nu>0$ for $t\in\left\langle t_B,t_E\right\rangle $  and  $y,w$ bounded is a positive function such that $H^\epsilon(t,y,w)$  is continuous; $\Delta(t)\geq\Delta>0$ on $\left\langle t_B,t_E\right\rangle$ where $\Delta$ is an arbitrarily small constant. These surfaces in $(t,y,w)$-space are defined by
\begin{equation*}
w=\pm\left( 2\left(H^\epsilon(t,y,w)-V(t,y)\right) \right)^{\frac12}
\end{equation*}
extending it as long as $w$ remains real. In our case such trajectories, lying on the surface $w=w(t,y,\epsilon),$ are bounded for every small $\epsilon$ (Fig. 2). On the charts for $K_1,$ and $K_3$ there is a motion in a single potential well and on the chart for $K_2,$ double well with a barrier in between.

Let $H^\epsilon(t)=H^\epsilon\left(t,y^\epsilon(t),w^\epsilon(t)\right),$  where $\left(y^\epsilon,w^\epsilon\right)$ is a solution of (\ref{Duff_osc_system1}), (\ref{Duff_osc_system2}) on $\left\langle t_B,t_E\right\rangle$ and let $y_L^\epsilon(t),$ $y_R^\epsilon(t)$ are the roots of equation 
\begin{equation*}
H^\epsilon(t)=V(t,y)
\end{equation*}
on $\left\langle t_B,t_E\right\rangle.$ 
Obviously,
{\setlength\arraycolsep{1pt}
\begin{eqnarray*}
&&y_L^\epsilon(t)<u_1(t)<y_R^\epsilon(t)\ \mathrm{on}\ K_1 \\
&&y_L^\epsilon(t)<u_1(t)<u_3(t)<y_R^\epsilon(t)\ \mathrm{on}\ K_2 \\
&&y_L^\epsilon(t)<u_3(t)<y_R^\epsilon(t)\ \mathrm{on}\ K_3 
\end{eqnarray*}}
Further, denote $y_L^0(t), y_R^0(t)$ the roots of equation
\begin{equation*}
H^0(t)+\Delta(t)=V(t,y).
\end{equation*}
Hence, $y_L^\epsilon(t)<u_2(t)<y_R^\epsilon(t)$ on $K_2$ and $y_L^\epsilon(t)\rightarrow y_L^0(t)$ from left side and $y_R^\epsilon(t)\rightarrow y_R^0(t)$ from right side on $K_1\cup K_2\cup K_3$ for $\epsilon\rightarrow 0^+.$

The derivative of $H^\epsilon(t)$ along any solution path of (\ref{Duff_osc_system1}), (\ref{Duff_osc_system2}), (\ref{Duff_osc_systemm3}) is
\begin{equation*}
H^{\epsilon'}(t)=w^\epsilon w^{\epsilon'}+f(y^\epsilon)y^{\epsilon'}-[m(t)y^\epsilon]'
\end{equation*}
\begin{equation*}
=w^\epsilon\left[-\frac{f(y^\epsilon)}{\epsilon a}+\frac{m(t)}{\epsilon a}-\frac{a'}{a}w^\epsilon\right]+f(y^\epsilon)y^{\epsilon'}-[m(t)y^\epsilon]'
\end{equation*}
\begin{equation*}
=-\frac{a'(t)}{a(t)}(w^{\epsilon})^2-m'(t)y^\epsilon.
\end{equation*}
For the regular $(\epsilon=1)$ and unfolded special case ($m(t)\equiv0$) the critical manifold is the union of the parallel straight lines on $\left(-\infty,\infty\right).$ Moreover, if $a$ is $C^1-$ function  with $a'(t)>0$ on $R,$ then the dynamics of dynamical system on  $K_2=\left(-\infty,\infty\right)$ (defining $t_{\min}=-\infty,$ $t_{\max}=\infty$) described by the equations (\ref{Duff_osc_system1}), (\ref{Duff_osc_system2}), (\ref{Duff_osc_systemm3}) is relatively simple, the $-\epsilon\frac{a'}{a}w$ term for increasing $a$ represents damping proportional to the velocity of particle. If started off with a certain amount of kinetic energy, the particle oscillates back and forth, gradually losing energy via damping and finally comes to rest at the bottom of one of the wells, for $t\rightarrow \infty$.
Further, for $a(t)\equiv1$ and $m(t)\equiv0,$ i.e. (\ref{Duff_osc_system1}), (\ref{Duff_osc_system2}), (\ref{Duff_osc_systemm3}) is conservative,  there exist the solutions switching between $y_L^\epsilon,$ $y_R^\epsilon$ on $K_2=\left(-\infty,\infty\right).$
\clearpage


\vglue2.5cm
\begin{figure}[ht]
\begin{center}
\includegraphics[scale=1, bb=6cm 16cm 15cm 20cm]{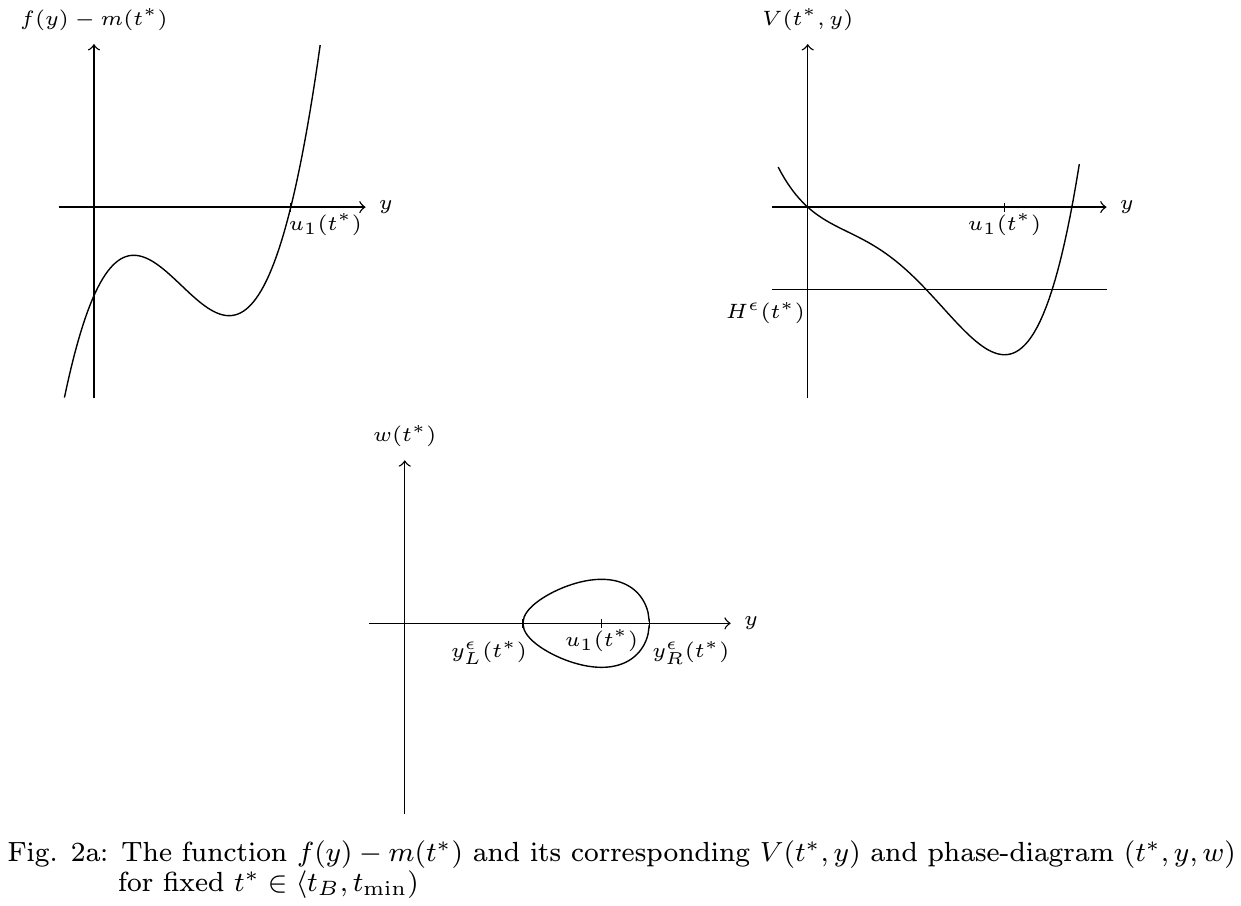}
\end{center}
\end{figure}

\vglue5.5cm
\begin{figure}[ht]
\begin{center}
\includegraphics[scale=1, bb=6cm 16cm 15cm 20cm]{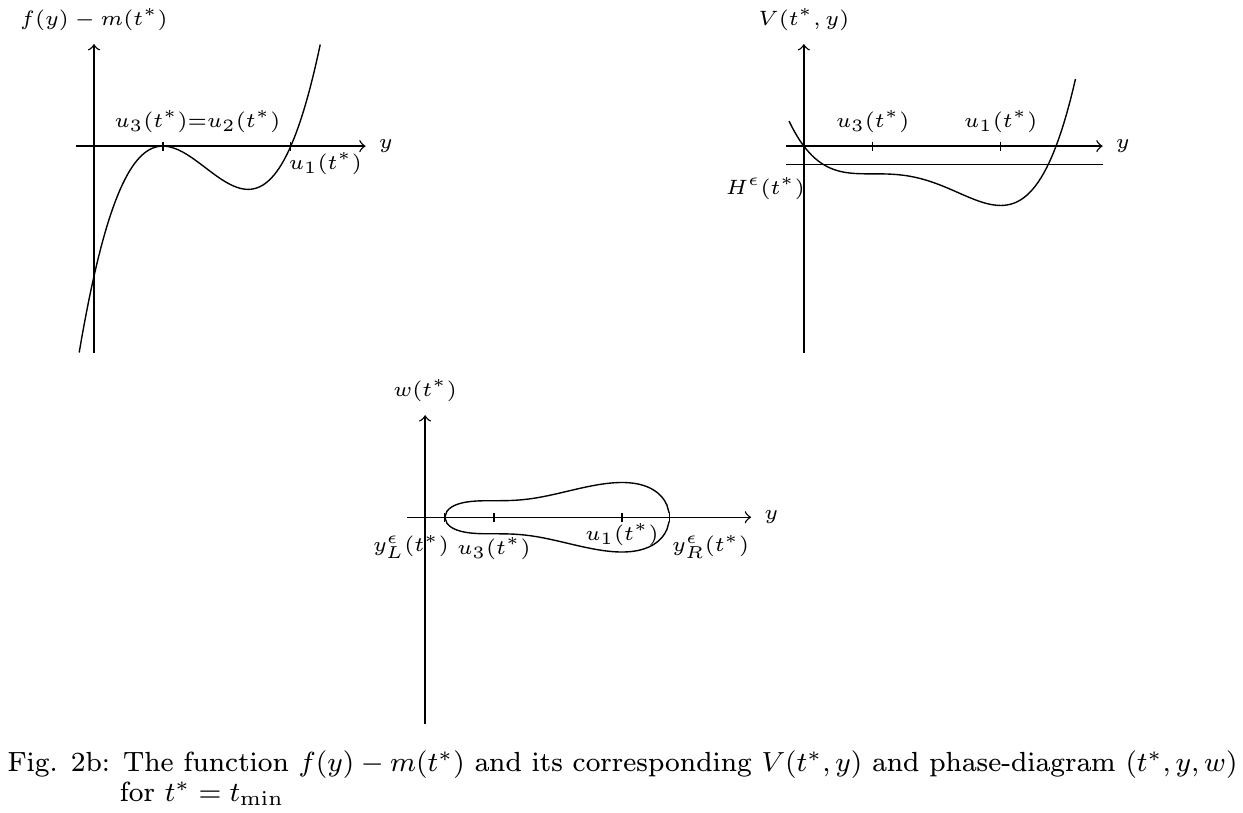}
\end{center}
\end{figure}

\clearpage

\vglue2.5cm
\begin{figure}[ht]
\begin{center}
\includegraphics[scale=1, bb=6cm 16cm 15cm 20cm]{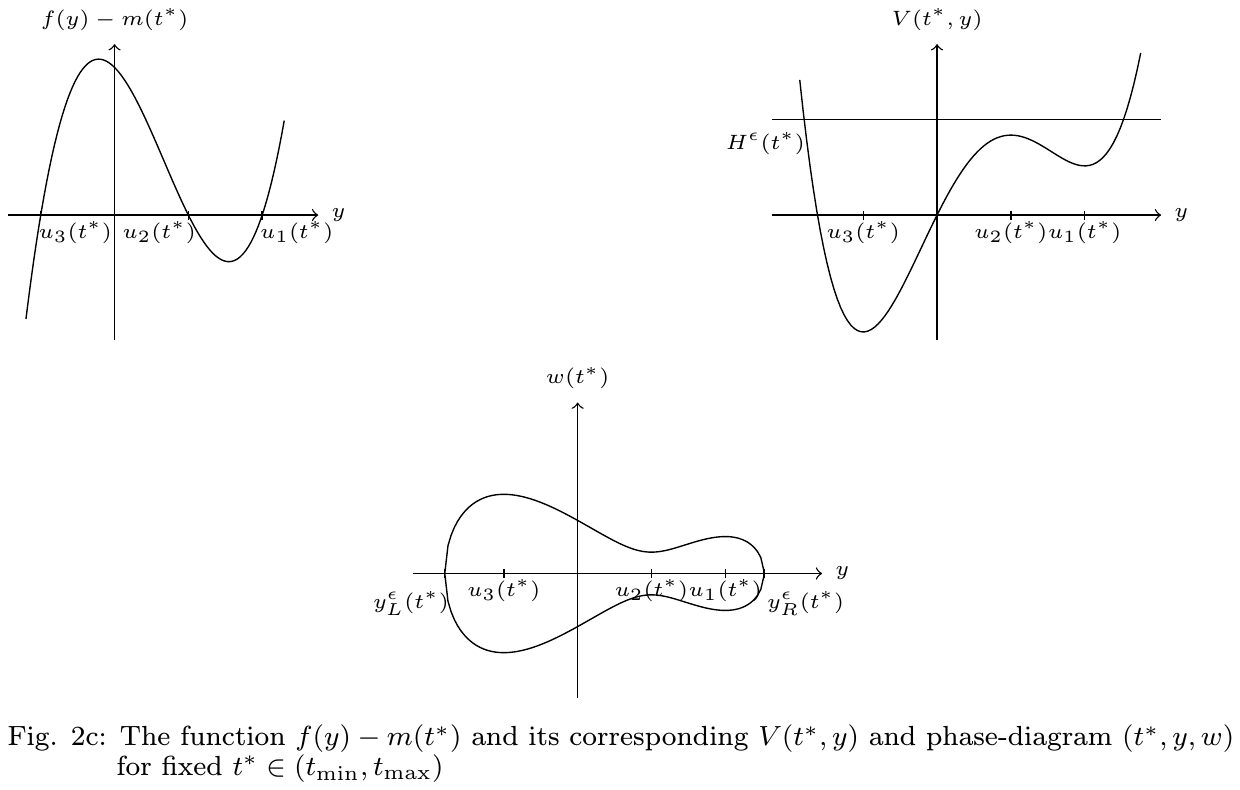}
\end{center}
\end{figure}
 
\vglue5.5cm
\begin{figure}[ht]
\begin{center}
\includegraphics[scale=1, bb=6cm 16cm 15cm 20cm]{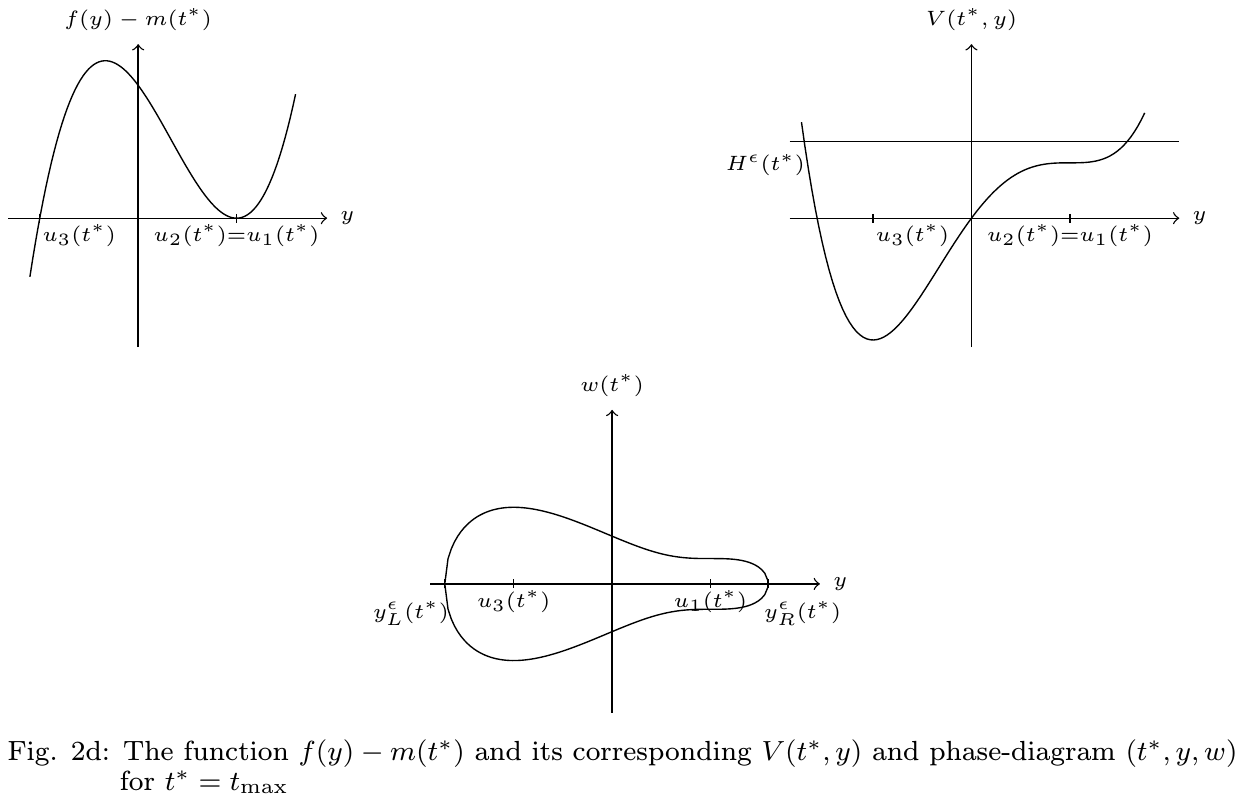} 
\end{center}
\end{figure} 
 
 \clearpage
 
\vglue2.5cm
\begin{figure}[ht]
\begin{center}
\includegraphics[scale=1, bb=6cm 18cm 15cm 20cm]{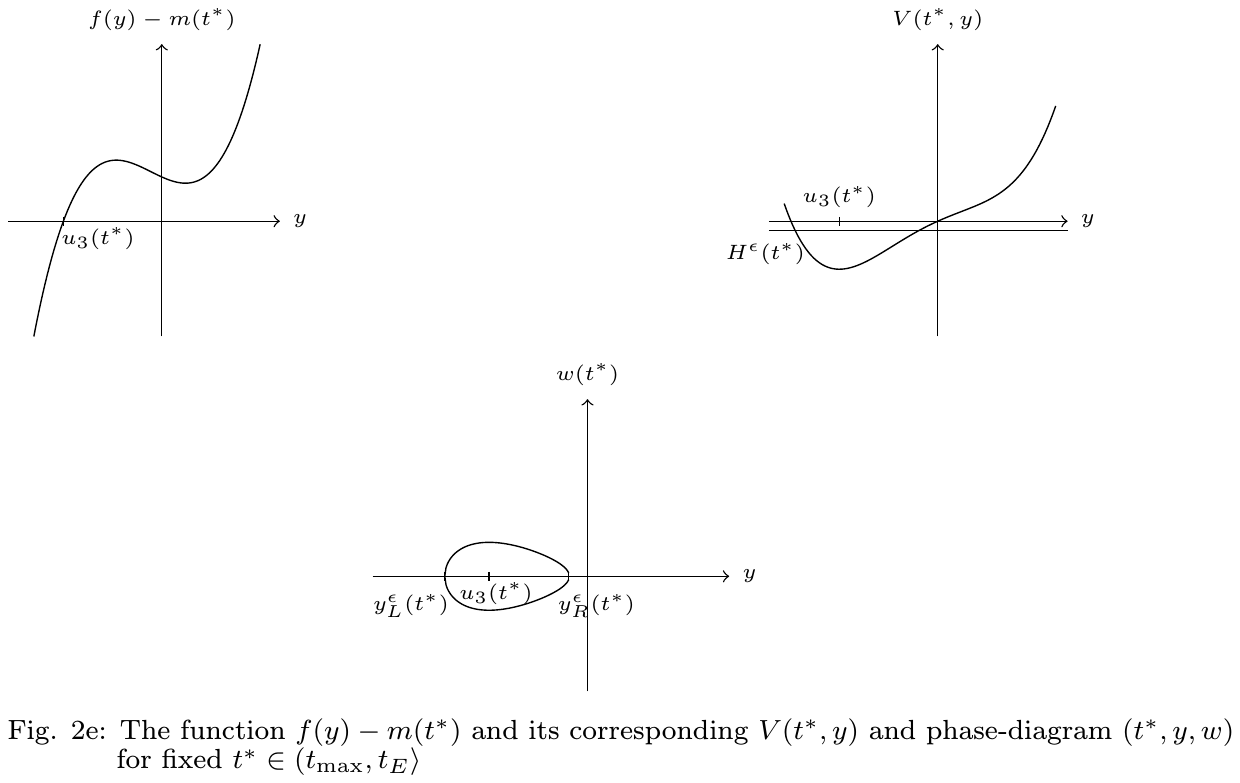} 
\end{center}
\end{figure} 
\vglue3.5cm
Let
\begin{equation*}\label{osc_criterion_f}
\chi(t,y)=\left(y-u_2(t)\right)\frac{\partial}{\partial y}\left[\frac{\int\limits_{u_2(t)}^y(f(s)-m(t))\mathrm{d}s}{\frac{\left(y-u_2(t)\right)^2}{2}}\right].
\end{equation*}
We make the following assumption.
\begin{itemize}
\item[$(A4)$]  The function
 \begin{equation*}\label{osc_criterion}
\chi(t,y)>-\frac{4\Delta}{\left(y-u_2(t)\right)^2}
\end{equation*}
for $y\in\left\langle u_1(t),u_2(t)\right)\cup\left(u_2(t),u_3(t)\right\rangle$ and for every fixed $t\in K_2.$ 
\end{itemize}

In this work we show that under the assumptions $(A1)-(A4)$ the Diff. Eq. (\ref{Duff_osc}) admits the nonlinear oscillations of solution $y^\epsilon$ with a frequency tending to infinity for $\epsilon\rightarrow 0^+.$
\section{Generalized polar coordinate transformation}
 We introduce the variable $v=\epsilon a^2y'$ and write (\ref{Duff_osc}) in the following system
{\setlength\arraycolsep{1pt}
\begin{eqnarray*}
& y'&=\frac v{\epsilon a^2}\label{Duff_osc_system3}\\  
& v'&=\frac{m(t)}{\epsilon }-\frac{f(y)}{\epsilon }.\label{Duff_osc_system4}
\end{eqnarray*}}
Then, we put $y=u_i(t)+r\cos\gamma,$ $i=1,2,3$ and $v=-r\sin\gamma$ on the charts $K_1,$ $K_2,$ $K_3,$ respectively.  We obtain the following differential equation for $\gamma$
\begin{equation*}\label{DE_gamma0}
\gamma'=\frac1\epsilon\left[\frac1{a^2(t)}\sin^2\gamma+\overline f_i(t,y)\cos^2\gamma+\frac{\epsilon u_i'(t)}{r}\sin\gamma\right]
\end{equation*}
or by using identity $\sin^2\alpha +\cos^2\alpha =1$
\begin{equation}\label{DE_gamma}
\gamma'=\frac1\epsilon\left[\frac1{a^2(t)}+\cos^2\gamma\left(\overline f_i(t,y)-\frac1{a^2(t)}\right)+\frac{\epsilon u_i'(t)}{r}\sin\gamma\right]
\end{equation}
where
\begin{equation*}\label{radius}
r=\sqrt{\left(y-u_i\right)^2+v^2}, 
\end{equation*}
\begin{equation*}
\overline f_i(t,y)=\frac{f(y)-m(t)}{y-u_i(t)},\quad \overline f_i\left(t,u_i(t)\right)=\frac{\mathrm{d}f}{\mathrm{d}y}\left(u_i(t)\right)\quad i=1,2,3.
\end{equation*}
\section{Analysis of chart for $K_2$ }
In this section we prove that under assumption $(A1)-(A4)$ is
\begin{equation}\label{gamma_ineq}
\gamma'\geq\frac1\epsilon c_{K_2}
\end{equation}
on $K_2,$ where $c_{K_2}$ is a positive constant.
First we estimate $r=r^\epsilon(t)$
{\setlength\arraycolsep{1pt}
\begin{eqnarray*}\label{radius2}
r_{\min}\left(K_i\right)&&\stackrel{\mathrm{def}}{=}\min\limits_{K_i}r^\epsilon(t) \\
&=&\min\limits_{K_i}\left\{u_i(t)-y_L^\epsilon(t), y_R^\epsilon(t)-u_i(t),\sqrt 2a(t)\sqrt{H^\epsilon(t)-V(t,u_i(t))}\right\}.
\end{eqnarray*}}
Because 
$$u_i(t)-y_L^\epsilon(t)>u_i(t)-y_L^0(t)>0,$$
$$y_R^\epsilon(t)-u_i(t)>y_R^0(t)-u_i(t)>0$$
and
$$H^\epsilon(t)-V(t,u_i)=\Delta(t)+h^\epsilon(t)>\Delta>0$$
is $r_{\min}\left(K_i\right)>0$ for every small $\epsilon$ and $i=1,2,3.$
Thus, third expression in (\ref{DE_gamma})
\begin{equation*}\label{ineq_2}
\left\vert\frac{\epsilon u_2'(t)}{r}\sin\gamma\right\vert\leq\frac{\epsilon\left\vert u_2'(t)\right\vert}{r_{\min}\left(K_2\right)}=O(\epsilon)
\end{equation*}
 on $K_2.$
For existence oscillations on chart $K_2$ is fundamental the analysis of the expression
\begin{equation}\label{expr1}
 \cos^2\gamma\left(\overline f_2(t,y)-\frac1{a^2(t)}\right).
\end{equation}
Clearly, 
\begin{equation*}\label{expr2}
\left\vert \cos^2\gamma\left(\overline f_2(t,y)-\frac1{a^2(t)}\right)\right\vert=\left\vert \frac{\left(y-u_2\right)^2\left(\overline f_2(t,y)-\frac1{a^2(t)}\right)}{\left(y-u_2\right)^2+a^2w^2}\right\vert.
\end{equation*}
For considered $H^\epsilon$ is $\left(y-u_2\right)^2+a^2w^2\neq0$ on the solution path therefore exists independent on $\epsilon$ constant $\delta_1=\delta_1(\eta)>0,$  such that for every $\eta,$ $0<\eta<\frac{1}{a^2(t)}$
\begin{equation*}
\left\vert \cos^2\gamma\left(\overline f_2(t,y)-\frac1{a^2(t)}\right)\right\vert\leq\frac1{a^2(t)}-\eta
\end{equation*}
for 
$y\in\left(u_2(t)-\delta_1,u_2(t)+\delta_1\right).$
Now analyze the expression  (\ref{expr1}) on the interval 
\begin{equation}\label{interval1}
\left(u_1(t)-\delta_2,u_2(t)-\delta_1\right)\cup \left(u_2(t)+\delta_1,u_3(t)+\delta_2\right),
\end{equation}
where $\delta_2$ is appropriate chosen positive constant such that $u_1(t)-\delta_2\geq y_L^\epsilon(t)$ and $u_3(t)+\delta_2\leq y_R^\epsilon(t).$

We obtain
\begin{equation*}\label{expr_ineq1}
-\cos^2\gamma\left(\overline f_2(t,y)-\frac1{a^2(t)}\right)=\frac{\left(y-u_2\right)^2\left(\frac1{a^2(t)}-\overline f_2(t,y)\right)}{\left(y-u_2\right)^2+2a^2(t)\left(H^\epsilon(t)-V(t,y)\right)}
\end{equation*}
\begin{equation*}
\leq\frac{\left(y-u_2\right)^2\left(\frac1{a^2(t)}-\overline f_2(t,y)\right)}{\left(y-u_2\right)^2+2a^2(t)\left(H^0(t)+\Delta-V(t,y)\right)}
\end{equation*}
\begin{equation*}\label{expr_ineq2}
=\frac1{a^2(t)}\left[\frac{\frac1{a^2(t)}-\overline f_2(t,y)}{\frac1{a^2(t)}-\frac{\int\limits_{u_2(t)}^y(f(s)-m(t))\mathrm{d}s-\Delta}{\frac{\left(y-u_2(t)\right)^2}{2}}}\right].
\end{equation*}
Because $\overline f_2(t,y)<0$ and $\int\limits_{u_2(t)}^y(f(s)-m(t))\mathrm{d}s<0$ on $\left(u_1(t),u_3(t)\right),$ $t\in K_2$ is the expression in  square brackets positive ($[\ ]>0,$ for short). Now we show that $[\ ]<1,$ independently on $\epsilon.$ From the assumption $(A4)$ we obtain that
\begin{equation*}\label{expr_ineq3}
\overline f_2(t,y)-\frac{\int\limits_{u_2(t)}^y(f(s)-m(t))\mathrm{d}s}{\frac{\left(y-u_2(t)\right)^2}{2}}>-\frac{\Delta}{\frac{\left(y-u_2(t)\right)^2}{2}}.
\end{equation*}
Hence
\begin{equation*}\label{expr_ineq4}
-\overline f_2(t,y)<-\frac{\int\limits_{u_2(t)}^y(f(s)-m(t))\mathrm{d}s-\Delta}{\frac{\left(y-u_2(t)\right)^2}{2}}.
\end{equation*}
Moreover, for $y=u_1(t), u_3(t)$ is $\overline f_2(t,u_1(t))\equiv\overline f_2(t,u_3(t))\equiv0$ on $K_2.$ Thus, there exists $\delta_2>0$ such that
\begin{equation*}\label{expr_ineq5}
-\frac{1}{a^2(t)}<\cos^2\gamma\left(\overline f_2(t,y)-\frac1{a^2(t)}\right)<0
\end{equation*}
on (\ref{interval1}).

Now, let $c_{K_2}$ from (\ref{gamma_ineq}) be $c_{K_2}=\min\left\{c_{K_2,1},c_{K_2,2}, c_{K_2,3}\right\}$ where
{\setlength\arraycolsep{1pt}
\begin{eqnarray*}
c_{K_2,1}=&&\min\left\{\eta-\frac{\epsilon\left\vert u_2'(t)\right\vert}{r_{\min}\left(K_2\right)},t\in K_2 \right\}\\
c_{K_2,2}=&&\min\left\{\frac{1}{a^2(t)}-\frac1{a^2(t)}\left[\frac{\frac1{a^2(t)}-\overline f_2(t,y)}{\frac1{a^2(t)}-\frac{\int\limits_{u_2(t)}^y(f(s)-m(t))\mathrm{d}s-\Delta}{\frac{\left(y-u_2(t)\right)^2}{2}}}\right]-\frac{\epsilon\left\vert u_2'(t)\right\vert}{r_{\min}\left(K_2\right)},\right.\\
&&\left. t\in K_2, y\in (\ref{interval1})\right\}            \\
c_{K_2,3}=&&\min\left\{\frac1{a^2(t)}\sin^2\gamma+\overline f_2(t,y)\cos^2\gamma-\frac{\epsilon\left\vert u_2'(t)\right\vert}{r_{\min}\left(K_2\right)},\right.\\
&&\left.t\in K_2, y\in \left\langle y_L^\epsilon(t),y_R^\epsilon(t)\right\rangle\setminus(\ref{interval1}),\gamma\in R\right\}.
\end{eqnarray*}}
Taking into consideration that $\overline f_2>0$ for $t\in K_2$ and $y\in \left\langle y_L^\epsilon(t),y_R^\epsilon(t)\right\rangle\setminus(\ref{interval1}),$ we conclude that $c_{K_2}>0$ for every sufficiently small $\epsilon,$ $\epsilon\in\left(\epsilon,\epsilon_0\right\rangle.$
\section{Analysis of the charts for $K_1$ and $K_3$ }
On the difference of $K_2,$ the analysis in the charts for $K_i, i=1,3$ is easy in comparation with $K_2$ one. The function $\overline f_i(t,y)>0$ for $t\in K_i$ and $y\in\left\langle y_L^\epsilon(t),y_R^\epsilon(t)\right\rangle,$ $i=1,3.$ Let
{\setlength\arraycolsep{1pt}
\begin{eqnarray*}\label{ck_K1K3}
c_{K_i}=&&\min\left\{\frac1{a^2(t)}\sin^2\gamma+\overline f_i(t,y)\cos^2\gamma-\frac{\epsilon\left\vert u_i'(t)\right\vert}{r_{\min}\left(K_i\right)},\right.\\
&&\left.t\in K_i, y\in \left\langle y_L^\epsilon(t),y_R^\epsilon(t)\right\rangle,\gamma\in R\right\},\ i=1,3.
\end{eqnarray*}}

The constants $c_{K_i},$ $i=1,3$ are positive for every sufficiently small $\epsilon,$ $\epsilon\in\left(\epsilon,\epsilon_0\right\rangle.$
Thus,  $\gamma=\gamma^\epsilon(t)$ is increasing on $K_i,$ $i=1,2,3$
\begin{equation}\label{gamma_ineqgen}
\gamma'\geq\frac1\epsilon c_{K_i}
\end{equation}

\section{Frequency control of  nonlinear oscillations}
In this section we show that the parameter $\epsilon$ play role modeling tool for the frequency control of the nonlinear oscillations.
Let us denote by $s_i$ the spacing between two succesive zeros of $y-u_i$ and by $z_i(y)$ the number of zeros of $y-u_i$ on $K_i,$ $i=1,2,3,$  where $y=y^\epsilon(t)$ is a solution of (\ref{Duff_osc}), then integrating the inequality  (\ref{gamma_ineqgen}) with respect to the variable $t$ between two succesive zeros of $y-u_i$  we obtain immediately
\begin{equation*}
\int\limits_{\mathrm{zero\ }(j)}^{\mathrm{zero\ }(j+1)}\gamma'\mathrm{d}t\geq\int\limits_{\mathrm{zero\ }(j)}^{\mathrm{zero\ }(j+1)}\frac{c_{K_i}}{\epsilon}\mathrm{d}t 
\end{equation*}
\begin{equation*}
\pi\geq\frac{c_{K_i}}{\epsilon}s_i.
\end{equation*}
Hence,
\begin{equation*}
s_i\leq\epsilon\frac{\pi}{c_{K_i}}
\end{equation*}
and
\begin{equation*}
\lim\limits_{\epsilon\rightarrow 0^+}z_i\left(y^\epsilon\right)=\infty,\ i=1,2,3.
\end{equation*}

Now we summarize the results of this article (pictorially, see Fig. 3).

\section{Statement of main result}

\begin{thm} \label{thm1}
Under the assumptions $(A1)$--$(A4)$ there exists solution $y^\epsilon$ of (\ref{Duff_osc}) for $\epsilon\in\left(0,\epsilon_0\right\rangle$ such that
$z_i\left(y^\epsilon\right)\rightarrow \infty$ with amplitude $y_R^\epsilon(t)-u_i(t)$ tendings to $y_R^0(t)-u_i(t)$ for subintervals of $K_i$  where $y^\epsilon-u_i\geq0$ and with amplitude $u_i(t)-y_L^\epsilon(t)$ tendings to $u_i(t)-y_L^0(t)$ for subintervals of $K_i$ where $y^\epsilon-u_i\leq0,$ $i=1,2,3.$
\end{thm}
\clearpage
\vglue3cm
\begin{figure}[ht]
\begin{center}
\includegraphics[scale=1, bb=6cm 17cm 15cm 20cm]{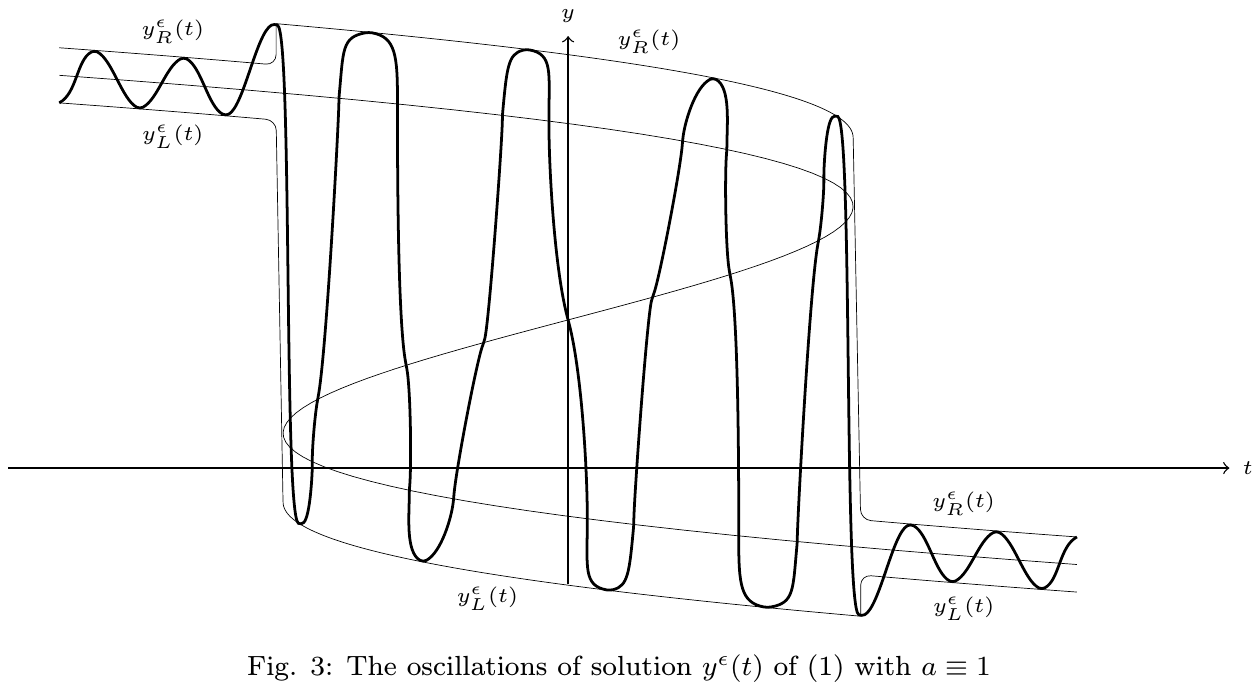}
\end{center}
\end{figure}

\vglue2cm

Finally, we remark that the proposed technique is an appropriate tool for detection and detailed analysis of the nonlinear oscillations in the dynamical systems but there is another powerful way to analyse the systems under consideration.   Indeed, after selecting the
new time $\tau=t/\epsilon$, system (\ref{Duff_osc_system1}), (\ref{Duff_osc_system2}), (\ref{Duff_osc_systemm3}) becomes a particular case of more general system of type
\begin{equation}\label{end_system1}
\frac{{\rm d}x}{\rm d\tau}=\frac{\partial H}{\partial y}(t,x,y)+\epsilon f_1(t,x,y),\quad 	
\frac{{\rm d}y}{\rm d\tau}=-\frac{\partial H}{\partial x}(t,x,y)+\epsilon f_2(t,x,y),\quad
t'=\epsilon.
\end{equation}
Assuming under the study of system (\ref{end_system1}) that for our values $t$ there exists a
family of closed trajectories inside the levels $\{(x, y) : H(t, x, y) = const\}$, one can
introduce new variables $(I,\phi)$ corresponding to these trajectories, in which the
subsystem
\begin{equation}\label{end_system2}
\frac{{\rm d}x}{\rm d\tau}=\frac{\partial H}{\partial y}(t,x,y),\quad 	
\frac{{\rm d}y}{\rm d\tau}=-\frac{\partial H}{\partial x}(t,x,y),\quad
t=const,
\end{equation}
takes the form
$$
\dot{I}=0,\quad \dot{\phi}=\omega(I,t),
$$
where $\omega(I,t)>0$. In new variables system (\ref{end_system1}) takes the form
\begin{equation}\label{end_system3}
\frac{{\rm d}I}{\rm d\tau}=\epsilon \Delta_1(I,\phi,t,\epsilon),\quad 	
\frac{\rm d\phi}{\rm d\tau}=\omega(I,t)+\epsilon \Delta_2(I,\phi,t,\epsilon),\quad
\frac{{\rm d}t}{\rm d\tau}=\epsilon.
\end{equation}
Now one needs to add only that to system (\ref{end_system3}) the standard averaging techniques
with respect $\phi$ could be applied (see e.g. \cite{GG}, \cite{KB}, \cite{SVM}).

\section*{Acknowledgments}
We would like to express our gratitude for all the valuable and constructive comments we
have received from the referee.

This research was supported by Slovak Grant Agency, Ministry of Education of Slovak Republic under grant number 1/0068/08.

\end{document}